\documentclass[12pt,a4paper]{amsart}

\usepackage{amssymb}
\usepackage{epsfig}

\begin{document}

\newtheorem{thrm}{Theorem}

\newtheorem{exa}{Example}

\newtheorem{thmf}{Th{\'e}or{\`e}me}

\newtheorem{thmm}{Theorem}

\def\thethmm{\ref{twist}$'$}

\newenvironment{thint}[1]{{\flushleft\sc{Th{\'e}or{\`e}me}}

      {#1}. \it}{\medskip} 

\newenvironment{thrm.}{{\flushleft\bf{Theorem}}. \it}{\medskip} 

\newenvironment{thrme}[1]{{\flushleft\sc{Theorem}}

      {#1}. \it}{\medskip} 

\newenvironment{propint}[1]{{\flushleft\sc{Proposition}}

      {#1}. \it}{\medskip} 

\newenvironment{corint}[1]{{\flushleft\sc{Corrolaire}}

      {#1}. \it}{\medskip} 

\newtheorem{cor}{Corollary}

\newtheorem{corf}{Corollaire}

\newtheorem{prop}{Proposition}

\newtheorem{defi}{Definition}

\newtheorem{rem}{Remark}

\newtheorem{deff}{D{\'e}finition}

\newtheorem{lem}{Lemma}

\newtheorem{lemf}{Lemme}

\newcommand{\sing}{\Sigma}

\newcommand{\zphi}{Z_{\phi}}

\newcommand{\new}{}

\def\fy{\varphi}

\def\ul{\underline}

\def\ol{\overline}

\def\obsf{{\flushleft\bf Remarque. }}

\def\obs{{\flushleft\bf Remark. }}

\def\tl{\widetilde}

\def\R{\mathbb{R}}

\def\comp{\mathbb{C}}

\def\C{\mathbb{C}}

\def\Z{\mathbb{Z}}

\def\HH{\mathbb{H}}

\def\N{\mathbb{N}}

\def\P{\mathbb{P}}

\def\Q{\mathbb{Q}}

\def\p{\pi_1(X)}

\def\Re{\mathrm{Re}}

\def\Im{\mathrm{Im}}

\def\H{\mathbf{H}}

\def\Hs{{Nil}^3}

\def\L{L_\alpha}

\def\h{\frac{1}{2}}

\def\M{\P(E)\times\P(E)^*\smallsetminus\mathcal{F}}

\def\m{\cp2\times{\cp2}^*\smallsetminus\mathcal{F}}

\def\g{\mathfrak{g}}

\def\l{\mathcal{L}}

\def\k{\mathfrak{h}}

\def\s3{\mathfrak{s}^3}

\def\nb{\nabla}

\def\Re{\mathrm{Re}}

\def\Im{\mathrm{Im}}

\newcommand\cp[1]{\mathbb{CP}^{#1}}

\renewcommand\o[1]{\mathcal{O}({#1})}

\renewcommand\d[1]{\partial_{#1}}

\def\sq{\square_X}

\def\gp{\dot\gamma}

\def\j{\mathcal{J}}

\def\sm{\smallsetminus}

\def\ra{\rightarrow}

\def\bi{\begin{enumerate}}

\def\ei{\end{enumerate}}

\newcommand{\pinfty}{p_{\infty}}

\def\so{\mathfrak{so}}

\def\1{\mathbf{-1}}

\title{A Singularity Theorem for Twistor Spinors} 

\author{Florin Belgun, Nicolas Ginoux, Hans-Bert Rademacher}

\thanks{AMS classification~: 53C21, 53A30, 32C10.\\
This work started while the second named author enjoyed the hospitality of the
Max-Planck Institute for Mathematics in the Sciences, which he would like to
thank for his support}
\date{2005--06--06}


\address{Florin Alexandru Belgun\\
Universit\"at Leipzig\\
Mathematisches Institut\\
Johannisgasse 26\\
D-04109 Leipzig}

 \email{Florin.Belgun@math.uni-leipzig.de}


\address{Nicolas Ginoux \\
Institut f\"ur Mathematik - Geometrie \\
Universit\"at Potsdam\\
Am Neuen Palais 10\\
D-14469 Potsdam}
\email{ginoux@math.uni-potsdam.de}


\address{Hans-Bert Rademacher \\
Universit\"at Leipzig\\
Mathematisches Institut\\
Johannisgasse 26\\
D-04109 Leipzig}
\email{Hans-Bert.Rademacher@math.uni-leipzig.de}

\begin{abstract}
We study spin structures on orbifolds. In particular, we show that if the
singular set has codimension greater than 2, an orbifold is spin if and only if its smooth part is.
On compact orbifolds, we show that any non-trivial twistor spinor admits at most one zero which is singular unless the orbifold is conformally equivalent to a round sphere. 
We show the sharpness of our results through examples. 
\vskip .6cm

\noindent 2000 {\it Mathematics Subject Classification}: Primary 53C27, Secondary 53C05,
53C24, 53C25, 14E15.

\medskip

\noindent{\it Keywords:} Orbifolds, Twistor spinors, ALE Spaces.
\end{abstract}

\maketitle

\section{Introduction}

\label{sec:introduction}


The twistor operator (also called the Penrose operator) and the Dirac operator
on a Riemannian spin manifold are obtained  by composing the Levi-Civita
covariant derivative with some natural linear maps. They are actually the two
{\it natural} first order linear differential operators on spinors. The
solutions of the corresponding P.D.E.'s (i.e., the kernels of these operators)
are the {\it twistor spinors} and, respectively,  the {\it harmonic spinors}, and they are both {\it conformally covariant}. Moreover, if we consider some appropriate {\it conformal weights}, they appear to be conformally invariant objects (as sections of some weighted spinor bundle).

\noindent However, it turns out that the {\it norm} of a twistor spinor defines
a special metric in the conformal class, for which the corresponding spinor is
actually parallel, or, more generally, sum of two {\it Killing spinors}. This
metric is then Einstein (of vanishing scalar curvature if and only if the
corresponding spinor is parallel).

\noindent Therefore, a dichotomy occurs: the study of the twistor spinors
without 
zeros reduces to the study of parallel resp. Killing spinors, 
and the study of twistor spinors with zeros which seems to be a 
``purely conformal'' problem. In both cases the corresponding 
manifolds and spinors can be described in full detail if we look for compact,
smooth, solutions.

\noindent The aim of this paper is to extend the above notions to the {\it
orbifold} case, 
where quotient-like singularities are allowed, and to concentrate on the 
twistor equation on compact orbifolds. Although the notions concerning spin
geometry 
on orbifolds may occur frequently as mathematical folklore, there is -- at our
knowledge --
very little written material focusing on the study of spinors on orbifolds 
(see, however, \cite{DLM} for a definition of spin orbifolds), and all questions
could be legitimate in this setting (in particular which topological
conditions
characterize a spin structure on an orbifold, see Theorem \ref{spinth}), 
however this type of quotient-like 
singularities naturally occur in the setting of twistor equations, as we explain below.

\noindent Indeed, using the solution of the Yamabe problem A. Lichnerowicz
proved 
that the only compact manifolds admitting twistor spinors with zeros are 
the standard spheres \cite{Li}, therefore conformally flat. Actually, 
the first example of a (non-compact) non conformally flat Riemannian spin
manifold 
carrying a twistor spinor with one zero was given
by W. K\"uhnel and H.B. Ra\-de\-ma\-cher \cite{KR95} using a
{\it conformal inversion} of the Eguchi-Hanson gravitational instanton
\cite{EH}. 

\noindent The resulting metric is invariant under an involution having a fixed
point, suggesting thus the possibility of extending the theory to orbifolds.

\noindent The main result of this paper, Theorem \ref{main}, establishes the
following facts:

\bi

\item On a compact orbifold (different from the standard sphere) admitting 
non-trivial solutions with non-empty vanishing set to the twistor equation, the
zero of such a spinor is 
unique and singular;
\item The (finite) cardinality of the {\it singularity group} of such a zero 
is not smaller than for any other point on the orbifold, and the equality case
characterizes the quotients of the standard sphere.

\ei

\noindent Furthermore, we give a series of examples of compact orbifolds 
carrying twistor spinors with zeros, cf. Section~\ref{sec:examples}. 
They rely on a general version of {\it conformal inversion} that allows 
us to compactify conformally an asymptotically locally Euclidean (ALE for short)
metric on some open manifold (or orbifold, in general) $M$ by adding one 
singular point at infinity, $p_\infty$. The resulting orbifold $N=M\cup
p_\infty$ 
admits then a twistor spinor with zero at $p_\infty$ if and only if $M$ admits a
parallel spinor, 
and therefore 3 possible holonomy groups occur for $M$: $SU(n)$ 
(in which case $M$ is a non-flat Ricci-flat ALE K\"ahler space, see below),
$Sp(n)$ (in which case $M$ should be ALE hyperk\"ahler). Such an example 
cannot be obtained by the below described construction of Joyce, 
\cite[Theorem 8.2.4]{Jo}), or when the holonomy group is $G_2$ (in which case
$M$ is 7-dimensional 
and thus $p_\infty$ is not an isolated singularity, therefore $M$ is not
smooth). 
In the present paper we give examples coming from the conformal 
compactification of a smooth manifold, therefore those examples belong all to
the first class.

\noindent Our fundamental examples can be then described as follows: 
We consider $S^{2n}$ as the conformal gluing of two copies of $\C^n$ (through
stereographical projections). 
The quotient of $S^{2n}$ by the standard action of a finite subgroup of $SU(n)$
is then an orbifold with 2 singular points, each of them being the origin in 
the corresponding chart. In one of the charts we resolve the singularity (we must assume there is a {\it crepant} resolution, see \cite[8.2]{Jo}) 
and get, in some cases\footnote{there is no criterium establishing that 
in general; that can be shown only on some explicit examples, cf. Section 
\ref{sec:examples}} a smooth open manifold $M$ 
admitting an ALE K\"ahler metric, see \cite{Jo}. Joyce's Theorem
\cite[Theorem 8.2.3]{Jo} 
states then the existence of a Ricci-flat ALE K\"ahler metric on $M$. 
Adding one point $p_\infty$ to $M$ or, equivalently, gluing to it the other 
copy of $\C^n/\Gamma$ gives us an orbifold $S^{2n}_\Gamma$. 
The resulting orbifold is thus obtained by resolving only one 
of the singular points of $S^{2n}/\Gamma$ (which obviously admits
itself two families of twistor spinors for the standard metric, each
of them vanishing in one of the two singularities).

\noindent We point out that, although the principle of conformal inversion holds
for 
all orbifolds in general, without restriction to the isolated singularity case,
there is no non-trivial example, so far, of Ricci-flat ALE structure on an 
orbifold with non-isolated singularities. In particular, the existence of 
odd-dimensional orbifolds admitting a twistor spinor with zero, as well as 
the existence of compact orbifolds admitting more than two linearly independent
such spinors remain open questions.\\

\noindent {\sc Acknowledgements } The authors would like to thank the anonymous
referee whose remarks led to an improvement of the material.

\section{Spin orbifolds}\label{sec:orbi}



In this section we will recall the basic properties of an orbifold and
define a spin structure on it.


\begin{defi} An $n$-dimensional {\em orbifold} $M$ is a Hausdorff
topological space, together with an atlas of charts $(U,f)$, where
$U$ is an open set in $M$ and 
$f:U\ra\R^n/\Gamma_f$ is a homeomorphism, where $\Gamma_f$
is a finite subgroup of $GL(n,\R)$ (acting on $\R^n$), such that the
transition functions $g\circ f^{-1}:f(U\cap V)\ra g(U\cap V)$ are
differentiable in the following sense: 
For any $x\in U\cap V$ there exists a small neighbourhood $W\subset U\cap V$ of
$x$
and a differentiable map (called the lift of the transition function)
$F:\pi_f^{-1}(f(W))\ra \pi_g^{-1}(g(W))$ such 
that $g\circ f^{-1}\circ \pi_f=\pi_g\circ F$, where $\pi_f:\R^n\ra\R^n/\Gamma_f$ and $\pi_g:\R^n\ra\R^n/\Gamma_g$ are the canonical 
projections.
\end{defi}

\obs Note that the lift $F$ is not unique unless $\pi_f^{-1}(f(W))$ and 
$\pi_g^{-1}(g(W))$ are connected, but it is always a local diffeomorphism.

\obs As any finite subgroup of $GL(n,\R)$ is conjugated to one sitting
in $O(n)$, we will suppose the groups $\Gamma$ above are
orthogonal. An {\em orientation} on $M$ is defined, as in the case of
a manifold, by the choice of an atlas of charts, such that the
Jacobian of the lift of transition functions has positive
determinant. This is well-defined if the groups $\Gamma$ actually lie
in $SO(n)$. As we are interested in oriented orbifolds, we will always
suppose the groups $\Gamma$ to lie in $SO(n)$.
\smallskip

\noindent The concept of a differentiable map between orbifolds can also be
defined as in the definition above. Of course, the
charts and their compositions with the canonical projections
$\pi_\Gamma:\R^n\ra\R^n/\Gamma$ are differentiable.






\begin{defi} The (singularity) group of a point $x\in M$ is the
conjugacy class of a minimal
(with respect to the inclusion) finite subgroup $\Gamma_x\subset SO(n)$ such that there exists a chart $(U,f)$ (which is then called {\em
minimal}) around $x$ with
$\Gamma_f=\Gamma_x$. If $\Gamma_x=\{1\}$ we say that $x$ is a smooth
point; otherwise it is called singular.
\end{defi}

\noindent If the singularity group of a singular point $x$  acts freely on 
$\R^n\sm\{0\}$, then the singularity is {\em isolated}, i.e., it is surrounded by smooth points. This can only happen in even dimensions (in the oriented case).

\obs A chart $(U,f)$ as in the definition above is called {\em minimal}
because there is no other chart around $x$ having a group $\Gamma'$
with less elements than $\Gamma_x$. Such a chart, composed with the
canonical projection from $\R^n$ to $\R^n/\Gamma$, yields a {\it ramified covering} of a neighbourhood of $x$ by an open set in
$\R^n$. A minimal chart yields, therefore, a $C^\infty$ map $\tilde f$ from a
neighbourhood of $0$ in $\R^n$ onto a neighborhood of $x$, such that
$\tilde f^{-1}(x)=\{0\}$.

{\noindent\bf Example. } In dimension 2, the group of a singularity is equal
to $\Z_n$, $n\ge 2$, and the ``total angle'' around such a point is
$2\pi/n$ (around a smooth point it is $2\pi$). So here the
singularities are always isolated (note that we restricted ourselves
to oriented orbifolds). A basic neighbourhood of such a ``conical''
point is actually homeomorphic to a disk, so every 2-dimensional
oriented orbifold is homeomorphic to a manifold. In larger dimensions 
the basic neighbourhood of an
isolated singularity is a cone over $S^{n-1}/\Gamma_x$, and, for
$n>2$, a quotient of a sphere is never homeomorphic to the sphere itself.

\medskip

\noindent If a singularity $x$ is not isolated, then the group of a neighbouring
point $y$ is the isotropy group of a point $z\in \tilde f^{-1}(y)\subset \R^n$ under
the action of $\Gamma_x$, so $\Gamma_y$ is isomorphic to a subgroup of 
$\Gamma_x$ (the minimal chart is obtained by restricting $\tilde f$ to a
$\Gamma_y$-invariant open subset $U$ of $z$, disjoint from
$\gamma(U)$, for any $\gamma\in \Gamma_x\sm \Gamma_y$). Then
$\Gamma_y$ will fix not only $z$, but actually a whole vector subspace
of $\R^n$, containing $z$. The set
of singularities of $M$ is then a (not necessarily disjoint) union of
orbifolds, each of even 
codimension (because of the orientability condition). 




\smallskip






{\noindent\bf Fundamental remark: } Any object on an orbifold $M$ can
be seen, in a neighbourhood of a point $x$, as a $\Gamma_x$-invariant
object on a local ramified covering by a smooth manifold (obtained
from a minimal chart composed with a canonical projection - we call
this {\it the (minimal) smooth covering of $M$ around $x$}).
\smallskip

\noindent We can now consider tensors on orbifolds: Locally they must come from
$\Gamma_x$-invariant tensors on the minimal smooth covering around
$x$.

\obs A vector field on an orbifold must vanish on any isolated
singularity and, in general, it must be tangent to the singular set,
because in these points $x$ it must be $\Gamma_x$-invariant. In
general, the tensor fields on $M$ must have particular
($\Gamma_x$-invariant) values in the
singular points. In particular, a {\it metric} on an orbifold is
locally a $\Gamma_x$-invariant metric on the minimal smooth covering
around $x$.

\noindent We can carry on most of the differential-geometric constructions on
orbifolds: for example, we can consider the Levi-Civita connection,
differentiate vector fields, take their Lie bracket etc. All those
operations may be performed locally in a chart, and there we will work
on the smooth covering with $\Gamma$-invariant objects.

\smallskip

\noindent We are interested now in putting a spin structure on an
orbifold. Before doing that, recall that the Spin bundle is a double
covering of the total space of the bundle of orthonormal frames, which
is non-trivial on each fiber. Note that, locally
around $x$, the frame ``bundle'' $SO(M)$ is just the quotient of the
frame bundle of the smooth covering under the action (by isometries)
of $\Gamma_x$. This action is always free, so the frame bundle of an
orbifold is a smooth manifold (but no longer a fiber bundle).

\noindent We will, however, continue to call $SO(M)$ the {\it bundle} of
orthonormal frames on $M$, and, in general, we will continue to use
the term {\it bundle} for quotients of (locally trivial) fiber bundles
on the local smooth coverings, and for objects which are locally of
this type.

\smallskip

\begin{defi} {\em (see also \cite{DLM})} A Spin structure on an orbifold $M$ is
given by a two-fold covering of the frame bundle $SO(M)$, which is non-trivial  over each fiber $SO_x$, $\forall x\in M$.
\end{defi}

\noindent We can then describe a spin structure on an orbifold as being locally a
$\Gamma_x$-invariant spin structure on the smooth covering around $x$,
but first we have to be able to lift the action of the group
$\Gamma_x$ of isometries to the Spin bundle of the smooth covering.

\begin{defi}\label{spinsing}
A singularity $x$ is said to be {\em spin} if there is a lift
$G_x\subset Spin(n)$ of $\Gamma_x\subset SO(n)$ which projects
isomorphically onto $\Gamma_x$ via the canonical projection from
$Spin(n)$ to $SO(n)$.
\end{defi}

\noindent If such a lift exists, it is not necessarily unique; for example, in
dimension 2, the group $\Z_n$ can be lifted to $Spin(2)$ if and only if
$n$ is odd. In that case there are two possible lifts. 

\obs There is another, deeper, motivation for the definition above: As
any spin structure on $M$ restricts to one on the smooth part $M\sm
S$, it should induce spin structures on the quotients of small spheres
around any singular point, i.e., on $S^{n-1}/\Gamma$. The condition
above is necessary and sufficient for  $S^{n-1}/\Gamma$ to be spin \cite[p.
47]{fr}.

\smallskip

\obs In \cite{fr} precise algebraic conditions on $\Gamma$ are
given for the existence of a lift $G$ of $\Gamma$ in $Spin(n)$. If
$n=4k+2$, $k\in \N^*$, such a lift is always unique, if it exists
\cite[sec. 2.2]{fr}.

{\noindent\bf Example} The group $\{\pm 1\}\subset SO(4)$ is the group
of a spin singularity: indeed, on $S^3/\{\pm 1\}\simeq \R P^3$ there are
exactly 2 (inequivalent) spin structures. On the other hand, there are 2
possible lifts of $\{\pm 1\}$ in $Spin(4)\simeq SU(2)\times SU(2)$,
which are $\{\pm 1\}\times \{1\}$ and $\{1\}\times \{\pm 1\}$.
\smallskip

\obs The action of $\Gamma$ by isometries on $SO(\R^n)$ commutes with
the right action of $SO(n)$ on the fibers of the frame bundle on
$\R^n$. The lifted action of $G$ on $Spin(\R^n)$ equally commutes
with the right action of $Spin(n)$ on the fibers. It follows that $G$ 
equally acts on every associated bundle $Spin(\R^n)\times_\rho F$,
where $\rho:Spin(n)\times F\ra F$ is a $C^\infty$ representation of
$Spin(n)$.

\noindent In particular, $G$ acts on the spinor bundles $\Sigma^\pm(\R^n)$ of
positive, resp. negative Weyl spinors of $\R^n$ (if $n$ is even; we
will mainly focus on this case).
That leads to the following definition (see also \cite{DLM}):

\begin{defi} The total spaces of the spinor bundles $\Sigma^\pm(M)$
 over an orbifold $M$ are the orbifolds obtained by gluing together the spinor bundles $\Sigma^\pm(M\sm S)$ with $\Sigma^\pm(U_x)\simeq\Sigma^\pm(\R^n)/G_x$, for all $x\in S$. A spinor field on an orbifold $M$ is
a pair of
smooth sections in each of these two bundles, such that the lifts in
any local smooth covering around  $x\in M$ are smooth ($G_x$-equivariant) spinors.
\end{defi}

\obs If the quotient  $\Sigma^\pm(\R^n)/G$ is not smooth, then the
value of any spinor field in $0$ must lie in the singular set, more
precisely in the set of fixed points of $G$ in $\Sigma^\pm(\R^n)$. We will
focus later on twistor spinors; they have the property that if the
spinor $\phi$ and the value of $D\phi$ (where $D$ is the Dirac
operator) simultaneously vanish in some point, then the twistor spinor
is everywhere zero (in a connected manifold), see next section. As
$\phi$ and $D\phi$ are sections in the 2 different spinor bundles,
they cannot both vanish at a singularity $x$ unless $\phi$ identically vanishes, which means first of all
that $G_x$ must have non-zero fixed points at least on $\Sigma^+(\R^n)$ or
on $\Sigma^-(\R^n)$. This implies certain constraints in dimension 4:

\begin{prop}
Let $\phi^+$ be a positive and $\phi^-$ a negative Weyl spinor field on a
4-dimensional orbifold. In a 
singular point $x$ at least one of them vanishes. Moreover,
they both vanish unless there is a complex structure on $\tl T_x M$
such that $\Gamma_x\subset SU(\tl T_xM)\simeq SU(2)$. In that latter
case, if $\Gamma_x\ne\{\pm 1\}$, then, for any other local spin
structure around $x$, every spinor field $\psi$ vanishes at $x$.
\end{prop}

\noindent The proof follows from the identity $Spin(4)\simeq SU(2)\times SU(2)$,
so any group $G\subset Spin(4)$ having a fixed point in $\Sigma^+$
must be totally contained in the second factor and conversely. If this
is the case, then any other lift of the corresponding projection
$\Gamma\subset SO(4)$ will act nontrivially on $\Sigma^+$ (always),
and on $\Sigma^-$ (unless $G=\Gamma=\{\pm 1\}$).
\endproof

















\section{A characterization of spin orbifolds}

The following result shows that is is enough to look at the smooth
part of an orbifold to see if it is spin and to determine its spin structure.

\begin{thrm}\label{spinth}
Let $M$ be an oriented orbifold, and let $S$ be the
set of its singularities. Assume $S$ is of codimension at least 4. Then $M$ is
spin if and only if the
manifold $M\sm S$ is spin. Moreover, the Spin structures on $M$ are in 1-1
correspondence with the spin structures on $M\sm S$. 
\end{thrm}

\noindent The codimension condition means the following: for any non-trivial 
element of any singularity group, its fixed point set is of codimension at least 4.
Note that this is essential: there are non-spin quotients of $\R^2$ by 
$\Gamma\subset SO(2)$, even if $\R^2\sm \{0\}/\Gamma$ is always a cylinder, 
therefore spin. On the other hand this codimension is always even for
oriented orbifolds.

\begin{proof}
The proof follows the following steps:

\bi
\item First we characterize topologically the total spaces of the frame bundles
of spin orbifolds (these are always smooth manifolds); we introduce the concept
of {\it orbifold universal covering}. In particular, extending a spin
structure
becomes equivalent to extending a certain discrete group action on a smooth
manifold.
\item We show that any spin structure on an $n$-sphere can be uniquely filled in to get one on the $n+1$-dimensional ball (the so-called {\em local model}).
\item We conclude using an induction on the dimension of the orbifold.
\ei

\noindent Recall that the frame bundle $SO(M)$ of an orbifold $M$, as defined before, is
smooth. 
Moreover, it comes with an action of the Lie algebra $\so(n)$, 
that induces a right action of the simply connected group $Spin(n)$ 
(recall that $n\ge 4$). This is actually a right $SO(n)$-action, therefore the
(non-trivial) element $\1$ lying in the kernel of the projection $Spin(n)\ra SO(n)$ 
acts trivially on $SO(M)$.

\begin{prop}\label{spin frame}
A spin structure on an orbifold $M$ is a double covering of $SO(M)$ on
which the element $\1\in Spin(n)$ acts freely. Moreover, $\1$ acts either
trivially (in which case $M$ is non-spin)
or freely on the universal covering $\tl{SO(M)}$, the action of $Spin(n)$ being
induced by the infinitesimal action of $\so(n)$. 
Therefore, a spin structure on $M$ is just a quotient of $\tl{SO(M)}$ by
a subgroup $G\subset \pi_1(SO(M))$, such that 
$$\pi_1(SO(M))=G\ltimes\{\mathbf{1},\1\}.$$
\end{prop}

\obs The advantage of this point of view is that it reduces the 
extension of a spin structure to the extension (by continuity, 
as we shall see) of the action of a group $G$ from an open dense set of
$\tl{SO(M)}$ to its whole. 

\begin{defi}
In both cases, the orbifold $\tilde{M}:=\tl{SO(M)}/Spin(n)$ 
is called the {\em orbifold universal covering} of $M$.
\end{defi}

\obs As its name suggests, the orbifold universal covering $\tl M$ can be
characterized by a universality property in the category of orbifolds,
similar with the one satisfied by the smooth universal covering of a
manifold (the coverings are to be understood here as orbifolds maps, 
so they are not necessarily locally invertible).

\noindent{\it Proof of Proposition}\ref{spin frame}
The claim is well-known if $M$ is smooth. For a general spin orbifold, 
we have defined a spin structure $Spin(M)$ to be a double-covering of $SO(M)$ 
which is fiberwise non-trivial. This implies that $\1$ has no fixed points 
over the smooth part of $M$. Now, a singularity $P$ of a spin orbifold is spin
itself,
and the restriction of the Spin-bundle $Spin(M)$ over a chart $U$ around $P$ is
a 
spin structure of an orbifold of the type $\R^n/\Gamma$, with
$\Gamma\subset SO(n)$. Let us describe now the spin structures of this
basic type of orbifold:

\begin{lem}\label{local} 
The spin structures on $B^n/\Gamma$ and on $S^{n-1}/\Gamma$ are both in 1-1 correspondence with the lifts $G$ of $\Gamma$ in $Spin(n)$. Moreover, any spin structure on $S^{n-1}/\Gamma$ can be uniquely
 filled in on $B^n/\Gamma$. Here, $B^n$ is the unit ball in $\R^n$, $\Gamma$ is the singularity group of $0$ and the codimension of the fixed point set of $\Gamma$ is at least $4$.
\end{lem}

\begin{proof} 
Let us first describe the $SO$-, resp. $Spin$-bundles on the orbifold 
$\R^n/\Gamma$. The orthogonal and spin frame bundles over $\R^n$, denoted by $SO(\R^n)$ and by $Spin(\R^n)$, respectively, are both Lie groups acting by isometries on $\R^n$ (actually $SO(\R^n)$ is the group of Euclidean transformations on $\R^n$), and the canonical projection from the last to the former is a group homomorphism.   
The group $\Gamma\subset SO(n)\subset SO(\R^n)$
is then a subgroup, and so is
$G\subset Spin(\R^n)$. As $SO(\R^n/\Gamma)=SO(\R^n)/\Gamma$ (here we
notice that $\Gamma$ acts on $SO(\R^n)$ by left multiplication --
hence commutes with the (right) $SO(n)$-action on the fibers), 
we see
that the fundamental group of  $M_{\Gamma}:=SO(\R^n/\Gamma)$ is the preimage
$\tl\Gamma\subset Spin(n)$ of $\Gamma$ under the fundamental
projection $p:Spin(n)\ra SO(n)$. So we have the following exact
sequence:
$$1\ra\{\pm1\}\ra\pi_1(M_{\Gamma})\ra\Gamma\ra 1.$$
On the other hand, if $\tl M_{\Gamma}$ is a two-fold covering of $M_{\Gamma}$,
it is
itself covered by $Spin(\R^n)$ which is the universal covering of both
$M_{\Gamma}$ and $\tl M_{\Gamma}$, so we have the exact sequence
$$1\ra G\ra\pi_1(M_{\Gamma})\ra\{\pm1\}\ra 1.$$
$\tl M_{\Gamma}$ is a spin structure on $\R^n/\Gamma$ if, moreover, the
two-fold covering $\tl M_{\Gamma}\ra M_{\Gamma}$ is non-trivial on each  fiber $SO_x$, $x\in\R^n/\Gamma$. This implies that $G$ cannot contain the
 non-trivial element in the $\Z_2$ from the first sequence, so the second sequence is a splitting of the first one.\\
\noindent So $G$ must be a lift of $\Gamma$. It is worth mentioning that
the double covering $\tl M_{\Gamma}$ of $M_{\Gamma}$ is actually determined by
the subgroup
$G$ of $\tl \Gamma$.
\smallskip

\noindent If we remove the singularity (and the corresponding fiber from the
frame bundles), the fundamental groups of the (smooth!) frame bundles 
remain the same, because
the codimension of the removed object is $n>2$. We can carry on the
same argument to conclude (see also 
\cite[Ch. 2.2]{fr} for the case when $\Gamma$ acts freely on $S^{n-1}$) that
there
is a 1-1 correspondence between the lifts $G$ of 
$\Gamma$ and the spin structures on $\left(\R^n\sm\{0\}\right)/\Gamma$, or,
equivalently, the spin structures on $S^{n-1}/\Gamma$. Note that the
$SO$, resp. $Spin$ bundle of $S^{n-1}$ can be identified with the Lie
group $SO(n)$, resp. $Spin(n)$.
\end{proof}

\noindent We first show now that, if $\1$ has some fixed point on $\tl{SO(M)}$,
then it must act trivially. We will show that $\tl{SO(M)}/\{\1,\mathbf{1}\}$ is
actually 
the frame bundle of $\tilde M$, therefore is smooth. Note that $\tl M$
has been defined, so far, as the space of orbits of $Spin(n)$ in
$\tl{SO(M)}$ and we can see, using the previous Lemma, that the
preimage $W_x$ in $\tl{SO(M)}$ of a small ball $U_x$ around $x\in M$ is
a union of copies of connected open sets, each of which is a covering
of $SO(U_x)$ and, therefore, diffeomorphic to a quotient of
$Spin(B_n)$ by a subgroup $H_x$ of $\tl\Gamma_x$, the $\Z_2$ extension of
the singularity group of $x$. If $\1\in H_x$ then it acts trivially on
$W_x$, if not, it has no fixed point there. So the fixed point set of
$\1$ is open and closed, therefore it is either empty or the total
space. In the second case, $M$ is clearly non-spin and $\tl SO(M)$ is
the frame bundle of $\tl M$; in the first case,
and again looking at the local model $W_x/Spin(n)$, which projects
over $U_x$, we can identify $\tl{SO(M)}/\{\1,\mathbf{1}\}$ with the frame bundle of $\tl M$, therefore the orbifold universal covering of $M$ is spin. 
Let $Spin(M)$ be a spin structure on $M$, i.e., a double covering of
 $SO(M)$ which is fiberwise non-trivial. Again, looking at the local
 model provided by the Lemma \ref{local}, we can see that $Spin(M)$ is
 locally a quotient of $Spin(B_n)$ by a subgroup of $\tl\Gamma$ of
 index 2, using the notations from Lemma \ref{local}. Again, if $\1$
 acts non-trivially (as it must be the case at least over the frame
 space of any smooth point), then it does not belong to this subgroup
 of index 2 and, therefore, acts freely on the quotient. So $\1$ has
 no fixed point on $Spin(M)$.
\endproof

\begin{lem}\label{univ}
If the codimension of the singular set of $M$ is at least $3$, then the
universal covering of $M\sm S$ is open dense in $\tl M$.
\end{lem}

\begin{proof} $\tl{SO(M)}$ is simply-connected and so is the preimage
  of $M\sm S$, because the preimage of $S$ is a union of submanifolds
  of codimension at least 3, as we can see from the local model, and
  removal of such submanifolds does not change the fundamental group
  of the total space.
But the preimage of $M\sm S$ is either the $SO$ or the $Spin$ bundle
of a smooth manifold open and dense in $\tl M$. This manifold is a 
covering of
$M\sm S$ and is also simply-connected, because its frame bundle
is. This proves the Lemma.
\end{proof}

\obs Note that $S$ is not necessarily a union of smooth submanifolds,
but the set of frames of all points in $S$ form a submanifold in
$SO(M)$. In $\tl M$, some of the singularities may be lifted (e.g., if
$M=\R^n/\Gamma$), so in
general $\tl{M\sm S}$ is just {\it contained} in the smooth part of
$\tl M$.
\smallskip

\noindent The extension of the spin structure on $M\sm S$ to $M$ proceeds now as follows: 
Let $Spin(M\sm S)$ be a spin structure on $M\sm S$. From Proposition 
\ref{spin frame} we know $Spin(M\sm S)=\tl{SO(M\sm S)}/G$, where 
$G\subset\pi_1(SO(M\sm S))$ and is isomorphic to $\Gamma=\pi_1(M\sm S)$. 
Because of the codimension condition, Lemma \ref{univ} 
states that $(\tl{SO(M\sm S)}$ is open dense in $\tl{SO(M)}=Spin(\tl{M})$, 
and actually $\Gamma$ acts on $\tl{M}$ by isometries, freely on 
the preimage of $(M\sm S)$ (which is therefore the universal covering of 
$M\sm S$).
We have to show that the action of $G$ can be uniquely extended by
continuity to
$Spin(\tl M)$. We will do this using Lemma \ref{local} and an induction 
on the dimension of $M$.
\smallskip

\noindent For $n=4$ the singular set $S$ is  a discrete set of points. Using the local
characterization of spin structures, as in Lemma \ref{local}, on
small balls around the points of $S$, we can extend these actions of
$\Gamma$, resp. $G$ on
the covering orbifold and on its spin bundle. This is just the
extension by continuity, but we can see it more precisely as follows:
Let $x$ be a singular point, and let $U_x$ be a small ball around
it. Then $U_x$ is covered by a union of copies of $V_x$, each of which
is itself diffeomorphic to a quotient of an Euclidean ball by a group
$\Gamma_0$ of isometries. Then we have 
$$\{1\}\ra\Gamma_0\ra\Gamma\ra\bar\Gamma\ra\{1\},$$
where $\bar\Gamma$ sends isometrically one copy of $V_x\sm\{x\}$ into
another, and $\Gamma_0$ preserves the copies of $V_x\sm\{x\}$. The
local model identifies $\Gamma_0$ with a subgroup of $SO(n)$ and the
extension of its action on whole $V_x$ is canonical. We can do the
same for $G$, and note that, in the local model, $G_0$ projects onto
$\Gamma_0$ and does not contain $\1$, actually $\1$ commutes with the
action of $G_0$ (before and after the extension), and acts freely on
the quotient $Spin(\tl M)/\Gamma$. We get, therefore, a spin structure on $M$
starting from one on $M\sm S$, as claimed.
We remark that the extension
procedure is the one used in Lemma \ref{local}, using here the fact
that the small geodesic spheres around singular points are smooth and
lie, therefore, in $M\sm S$.
\smallskip

\noindent For the induction step we suppose that on any $(n-1)$-dimensional
oriented orbifold whose singular set is of codimension at least 4, and
such that the smooth part of it is spin, the action of the discrete
group $G$ as above can be extended to the singular set.
\smallskip

\noindent We will use this induction hypothesis only for the small geodesic
spheres around singular points. Let $x\in S$, then we get an extension
of the action of $G$ on the preimage in $Spin(\tl M)$ of a small
sphere around $x$, because that sphere quotient is of dimension $(n-1)$ and the
action of $G$ is well-defined above the smooth part of it. But then we can fill
in the obtained spin structure on this sphere to get one on the
corresponding small ball around $x$, using again Lemma \ref{local}. By
doing so, we actually extend the action of $G$, by
continuity, over the preimage of $x$ in $Spin(\tl M)$. As this
extension is unique, it is independent of the choices 
made and we obtain a free action of $G$ on $Spin(\tl M)$, such that $\1$
acts freely on the quotient, and such that the further quotient under
$\1$ is $SO(M)$. Therefore, $Spin(\tl M)/G$ is a spin structure on $M$.
\end{proof}

\section{Zeros of twistor spinors and conformal inversion}\label{sec:conf}


Let us first recall the twistor equation on a Riemannian spin manifold
of dimension $n: $ 
We call a spinor field $\phi$ a {\em twistor spinor} if the
following {\em twistor equation} holds for all tangent vectors $X:$
\begin{equation}\label{eq:twistor}
\nabla_X\phi+\frac{1}{n} X\cdot D \phi=0\,.
\end{equation}

\noindent Here $\nabla$ is the {\em spin-connection} on the spinor bundle
$\Sigma (M),$ the dot ``$\,\cdot\,$'' denotes the Clifford multiplication and
$D$ the Dirac operator.
This definition extends to the orbifold case: around a singularity $x$ the
spinor field $\phi$ 
and tangent vectors $X$ must then have
$G_x$-equivariant liftings. The spin-connection $\nabla_X\phi$  and 
the Dirac operator $D \phi=\sum_{i=1}^n e_i \cdot \nabla_i \phi$ of
a $G_x$-equivariant spinor field $\phi$ are again $G_x$-equivariant.  
If $\phi$ is a twistor spinor then one computes for the 
derivative of the Dirac operator $D \phi:$

\begin{equation}\label{eq:deriv-tw}
\nabla_XD\phi= \frac{n}{2}L(X)\cdot \phi\,.
\end{equation}

\noindent Here $L$ is the $(1,1)-${\em Schouten tensor} defined by
$$ L(X)= \frac{1}{n-2}\left(\frac{s}{2(n-1)}X - {\rm Ric} (X)\right)$$
with the {\em Ricci tensor} ${\rm Ric}$ and the scalar
curvature $s.$ One can view Equation~(\ref{eq:twistor})
and Equation~(\ref{eq:deriv-tw}) as a parallelism condition:
We define on 
the double spinor bundle $E= \Sigma (M) \oplus \Sigma (M)$ 
the connection
$$ \nabla_X^E=\left(\begin{array}{cc}
\nabla_X & \frac{1}{n}X \cdot \\
-\frac{n}{2}L(X) \cdot & \nabla_X
\end{array}
\right)\,,$$
i.e., 
$$\nabla_X^E(\phi,\psi)=\left(\nabla_X\phi+\frac{1}{n}X\cdot \psi,
-\frac{n}{2}L(X) \cdot \phi + \nabla_X \psi\right) \,.$$
Then we obtain 

\begin{lem}\label{lem:tw-parallel}\cite[ch. 1.4, Thm. 4]{BFGK}.
A twistor spinor $\phi$ on a Riemannian spin manifold is 
uniquely determined by a parallel section of the bundle
$E$ with connection $\nabla^E.$ More precisely, if 
$\phi$ is a twistor spinor, then $(\phi, D\phi)$ is a parallel
section of $(E,\nabla^E)$ and if $(\phi,\psi)$ is a parallel
section of $(E,\nabla^E)$ then $\phi$ is a twistor spinor and 
$\psi=D\phi.$
\end{lem}

\begin{cor}\label{isol}
Any zero $P$ of a twistor spinor $\phi$ is isolated, 
more precisely $\nabla_X\phi(P)\ne 0$ for any non-zero vector $X\in T_PM$.
\end{cor}

\begin{proof}
If both $\phi$ and $D\phi$ vanish at $P$, then $\phi$ is identically zero
(as unique parallel section in $E$ with this initial data). Therefore
$\nabla_X\psi=-\frac{1}{n} X\cdot D\phi$ is non-zero if $X\ne 0$, as Clifford product of a 
non-zero vector 
with a non-zero spinor $D\phi$.
\end{proof}


\noindent We use the following notation for the open Euclidean ball
of radius $R:$
$B_R:=\{x\in \R^n\,|\,\|x\|< R\},$ then
$\overline{B}_R:=\{x \in \R^n\,|\, \|x\|\le R\}$ is the 
corresponding closed
Euclidean ball.

\begin{defi}\label{def:ale}
{\em (a). }
Let $\Gamma$ be a finite subgroup of $SO(n)$ that acts
freely on $\R^n\sm\{0\}$ and let $U$ be an open subset of 
a Riemannian manifold. Then the open subset $U \subset M$ carries an
{\em Asymptotically Locally Euclidean (ALE) coordinate system}
$y=(y_1,\ldots,y_n)$
of order $(\tau, \mu)$ if the following conditions 
are satisfied:

There is for some $R>0$ a diffeomorphism
$y \in \left(\R^n-\overline{B}_R\right)/\Gamma\mapsto \phi (y) \in U$
such that the metric coefficients
$g_{ij}(y)=g\left(\frac{\partial}{\partial y_i},
\frac{\partial}{\partial y_j}\right)$ 
with respect to the coordinates $y=(y_1,\ldots,y_n)$
and its derivatives $$\partial_{i_1 \ldots i_k}g_{ij}=
\frac{\partial^k}{\partial y_{i_1}\cdots
\partial y_{i_k}}g_{ij}$$ 
have the following asymptotic behaviour for 
$\rho=\|y\|=\sqrt{\sum_1^n y_i^2} \to \infty :$
$$ g_{ij}-\delta_{ij}=O(\rho^{-\tau}) ; 
\partial_{i_1\ldots i_k}g_{ij}=O\left(\rho^{-\tau-k}\right)$$
for all $k=1,2,\ldots,\mu.$ If the
group $\Gamma=\{1\}$ is trivial, then the coordinate system is
called {\em Asymptotically Euclidean.}


\smallskip

{\noindent{\em (b). }}
\noindent We call a non-compact Riemannian manifold $\overline M$ 
of dimension $n$ 
{\em Asymptotically Locally Euclidean} or short
{\em ALE} of order $(\tau, \mu)$ if 
there is a compact subset $M_0$ such that the complement
$M-M_0$ carries an asymptotically Euclidean coordinate system
of order $(\tau, \mu).$
\end{defi}






\obs There is a close 
relationship between manifolds carrying twistor spinors
with zeros and non-compact manifolds with parallel spinors
with an end carrying an ALE coordinate system, more precisely:


\noindent Let $(M,g)$ be a Riemannian
spin manifold with a twistor spinor $\phi$ having a zero $p$. 
It is isolated and the length $\|\phi\|$ behaves like a distance function
in the neighbourhood of $p$ (see 
Corollary \ref{isol}). 
Given normal coordinates 
$x=(x_1,\ldots,x_n)\in B_{\epsilon}\subset\R^n$ 
in  a neighbourhood $U$ of $p$
we define {\em inverted normal coordinates}
$y=(y_1,\ldots,y_n)\in \R^n-B_{\epsilon^{-1}}, y=x/\|x\|^2.$
Then the conformally equivalent metric
$\left(U-\{p\},\overline{g}=g/\|\phi\|^4\right)$
carries a parallel spinor and an asymptotically
Euclidean coordinate system of order $(3,2).$
It is locally irreducible unless it is flat (see
\cite[Theorem 1.2]{KR98} for details).



\noindent On the other hand one can use a metric with parallel spinors
having an end with an ALE coordinate system to produce examples
of twistor spinors with zeros:


\begin{lem}\label{lem:metricconf}
If on the open subset $U$ diffeomorphic
to $\left(\R^n\sm\overline{B}_R\right)/\Gamma$ of a smooth, i.e. 
$C^{\infty}$ Riemannian manifold
$(M,g)$ there is an ALE-coordinate system $y$ with radius function
$\rho=\|y\|=\sqrt{\sum_1^ny_i^2}$ of order $(\tau,\mu)$
with $\mu \ge \tau-1 \ge 2$,
then the conformally equivalent metric 
$\overline{g}=\rho^{-4}g$ extends as a $C^{\tau-1}$ metric
to the one-point
completion $U \cup {\pinfty}$ diffeomorphic
to $B_R/\Gamma.$ 
\end{lem}


\begin{proof}
We denote by $y=(y_1,\ldots, y_n),\rho >R$ 
asymptotically Euclidean coordinates
and denote
$$g_{ij}(y)=g\left(\frac{\partial}{\partial y_i},
\frac{\partial}{\partial y_j}\right)=\delta_{ij}+h_{ij}(y).$$
Then $h_{ij}=O(\rho^{-\tau}),
\frac{\partial^k}{\partial y_{i_1}\cdots y_{i_k}}h_{ij}(y)=
O(\rho^{-\tau-k})$ for all $k=1,2,\ldots,\mu.$
We use the {\em inversion} $z=\rho^{-2}y$ and obtain
with the formula
$\frac{\partial}{\partial z_i}=\rho^{-2}\frac{\partial}{\partial y_i}
-2 \frac{z_i}{\rho^2}\sum_{k=1}^n z_k \frac{\partial}{\partial y_k} $
for the coefficients $\overline{g}_{ij}$
of the
conformally equivalent metric $\overline{g}$
with respect to the inverted coordinates $z=(z_1,\ldots,z_n);$
$z_i=\rho^{-2}y_i:$

\begin{eqnarray}
\overline{g}_{ij}(z)&=&\hat g \left(\frac{\partial}{\partial z_i},
\frac{\partial}{\partial z_j}\right)=
\delta_{ij}+h_{ij}-\frac{2}{\rho^2}\left(
z_i \sum_k z_k h_{kj}+z_j\sum_l z_lh_{il}\right)\nonumber \\
&& + \frac{4}{\rho^4} z_iz_j \sum_{k,l} z_kz_l h_{kl} \,.
\end{eqnarray}

\noindent It follows that $\overline{g}_{ij}(z)=\delta_{ij}+O(r^{\tau})$
with $r=\|z\|=\sqrt{\sum_1^n z_i^2}=\rho^{-1}$
and 
$\frac{\partial^k}{\partial z_{i_1}\cdots\partial z_{i_k}}
\overline{g_{ij}(z)}
= O\left(r^{\tau - k}\right).$

\noindent Hence we obtain that the function $\overline{g}_{ij}(z)$
extends to $z=0$ as a $(\tau -1)$-times 
continuously differentiable function.
\end{proof}


\noindent The following result establishes a {\it conformal completion} of a Ricci-flat ALE K\"ahler metric which will provide us with 
examples of conformal orbifolds admitting twistor spinors with zero.

\begin{thrm}\label{thm:parconf}
Let $\Gamma$ be a subgroup of $SO(n), n=2m$ acting
freely on $\R^n\sm\{0\}$ and let $(M,g)$ be an ALE
Riemannian spin manifold $(M,g)$ 
(with a $C^{\infty}$-metric) of order $(\tau, \tau)$ (i.e. asymptotically
Euclidean to $\R^n/\Gamma$ at infinity) with $\tau \ge 2$
and holonomy group ${\rm SU (m)}.$
Then there is a one point 
conformal completion 
$N=M \cup \{\pinfty\}$ of $(M,g)$ to a compact
Riemannian spin orbifold with singular point $\pinfty$
whose singularity group is $\Gamma.$
The metric $\overline{g}$ on $M=N-\{\pinfty\}$ is 
$C^{\infty}$-smooth, conformally equivalent to
$g$ on $M=N-\{\pinfty\}$ and it is a $C^{\tau-1}$ metric on $N.$\\
Then there is a spin structure on the orbifold $N$
with a two-dimensional space of twistor spinors, and 
all nontrivial twistor spinors have exactly one zero point,
which is the singularity point $\pinfty.$
\end{thrm} 


\begin{proof}
Since the holonomy group is ${\rm SU}(m)$ one can conclude that
the manifold $M$ is spin and has a preferred spin structure
for which the space of parallel spinors is two-dimensional,
cf. \cite[Corollary 3.6.3]{Jo}.
It follows from Lemma~\ref{lem:metricconf} 
that the conformally equivalent
metric $\rho^{-4}g$ can be extended to the orbifold 
$N=M \cup \{\pinfty\}.$ 
If $U$ is a sufficiently small neighbourhood of $\pinfty$ then
there is a Riemannian covering 
$(V\sm\{p\},\overline{g})$
of $(U\sm\{p\},\overline{g})$ 
with Riemannian covering group $\Gamma$ which is diffeomorphic to
$B^n_1\sm\{0\}$
carrying a twistor spinor $\phi$ invariant under $\Gamma.$
It follows from the behaviour of twistor spinors under conformal
changes that its norm with respect to the Hermitian
metric on the spinor bundle of $\left(V\sm\{p\},\overline{g}\right)$
is given by $\|\phi\|=r,$ i.e. the spinor field $\phi$
can be continuously extended to $U$ by setting
$\phi (\pinfty)=0.$ 
As pointed out above a twistor spinor $\phi$ is in one-to-one
correspondence with a parallel section of the bundle
$E$ with connection $\nabla^E.$ Then it follows from the next
Lemma~\ref{lem:extension} 
that it extends to a continuously differentiable section,
hence it is a $\Gamma$-invariant twistor spinor on 
$(V,\overline{g}).$ Therefore the orbifold $(N,\overline{g})$ 
carries a twistor spinor with zero in $p.$
\end{proof}


\noindent In the proof of Lemma~\ref{lem:metricconf} we used the following
general result to extend the twistor spinor into the singularity:


\begin{lem}
\label{lem:extension}
Let $E\ra M$ be a $\mathcal{C}^1$ vector bundle equipped with a
continuous linear connection (i.e., in a $\mathcal{C}^1$ trivialization map,
the coefficients of the connection form are continuous forms on
$M$, with values in the set of linear endomorphisms of the fiber), and
let $\sigma$ be a parallel section over $M\sm\{p\}$. Then $\sigma$ can
be uniquely extended to a $\mathcal{C}^1$ global section (hence parallel).
\end{lem}


\begin{proof}
By replacing $M$, if necessary, with a smaller neighbourhood of $p$, we
can suppose the bundle is trivial, so $E=M\times F$, and the
connection form is a continuous 1-form $\omega$ on $M$, with
values in $End(F)$. Consider some Riemannian metric on $M$ and an
Euclidean metric on $F$. Without loss of generality we can assume that
$\omega_p=0$, and (by restricting it, if necessary, to a smaller
neighbourhood around $p$),
\begin{equation}\label{norm}
|\omega_x(X)\cdot V|\le \varepsilon |X| |V|,
\end{equation}
where the dot means matrix multiplication, and the norms of $X\in
T_xM$, $V\in F$ are computed using the above chosen Riemannian,
resp. Euclidean Metric on $M$ and $F$.
We consider a parallel lift of a $\mathcal{C}^1$ loop $c$ in $M$,
passing through $x$, and contained in a compact ball around $x$. It is
a curve $t\mapsto (c(t),\gamma(t))$ in $M\times F$. This
lift is characterized by the following linear ODE:
$$\gamma'(t)=-\omega_{c(t)}(c'(t))\cdot\gamma(t)$$
Suppose $|c'|\equiv 1$, so $c$ is an arc length parameterization. By
taking the scalar product in $F$ with $\gamma(t)$, and using
(\ref{norm}), we obtain
$$|(|\gamma(t)|^2)'|\le 2\varepsilon |\gamma(t)|^2,$$
hence 
$$|\gamma(t)|\le e^{\varepsilon t}|\gamma(0)|,$$
and 
\begin{equation}\label{gam}
|\gamma'(t)|\le \varepsilon e^{\varepsilon t}|\gamma(0)|.
\end{equation}
This gives us a bound for the length of $\gamma$, depending
continuously on the length of $c$, and linearly on the initial data
$\gamma(0)$. The first conclusion is that the lift is defined for all
times.
Let $\sigma:M\sm\{p\}\ra F$ be a parallel section outside $p$. Let $c$
be an arbitrary $\mathcal{C}^1$ curve on $M$, such that $c(0)=p$ and
$c'(0)\ne 0$. Starting from $c(t_0)\ne p$, for some small time $t_0$,
we consider the lift $t\mapsto (c(t),\sigma(c(t)))$, which is
horizontal, because $\sigma$ is parallel. But we can define this lift
for all times, including $0$. We obtain a curve $\gamma$, depending on
$c$, and a point $(p,\gamma(0))=(p,A)$ in the fiber over $p$.
We want to show that this point does not depend on the choice of the
curve $c$. Suppose that, for another curve $\tilde c$, equally arc length
parametrized, we obtain a lift $\tilde\gamma$, such that $\tilde\gamma(0)=B\ne
A$.
Using $c$ for the beginning, $-\tilde c$ for the end (i.e., we run
$\tilde c$ in
the reverse sense), and connecting them smoothly, we can get, for any
$n\in\N^*$, a sequence of $\mathcal{C}^1$ arc length parametrized
loops $c_n$ in $M$, of length $l_n\le 1/n$, such that
$c_n(0)=c_n(l_n)=p$.
The corresponding lifts $\gamma_n$ through $\sigma$ are well defined over
$M\sm\{p\}$, and because of the argument above,
$\gamma_n(0)=\gamma(0)=A$ (the curves $c$ and $c(n)$ coincide around
$0$), and $\gamma_n(l_n)=\tilde\gamma(0)$ (the curves $-\tilde c$ and
$c_n$ coincide around $0$, resp. $l_n$).
From (\ref{gam}), the length of $\gamma_n$ is smaller than $|A|
(e^{\varepsilon/n}-1)$, therefore the distance between $A$ and $B$ must be smaller that this expression.\\
So $A=B$.\\
This allows us to define $\sigma(p):=A$ and get a continuous global
section. Because the lifts through $\sigma$ of $\mathcal{C}^1$ curves
through $p$ are parallel lifts through $A$, we easily obtain the
differentiability $\sigma$ in $p$, and actually that $\sigma$ is
$\mathcal{C}^1$.
It is also, by construction, parallel.
\end{proof}

\section{Main results}\label{sec:main}

\begin{thrm}\label{main}
Let $M$ be a compact Riemannian spin orbifold. If $M$ carries a non-trivial
twistor spinor with zero
at $p_0$ then this zero is unique
and $p_0$ is a singular point unless the orbifold is conformally 
equivalent to a 
round sphere. In addition, for every point $q \not=p_0$
the order $\# \Gamma_q$ of the
singularity group satisfies $$ \# \Gamma_q \le \# \Gamma_{p_0}$$
with equality only if the orbifold is a quotient of a standard sphere,
hence conformally flat.
\end{thrm}

\noindent Then the conformally flat case is discussed, under some assumptions, in Proposition \ref{confflatex}.

\begin{proof}
Let $\phi$ be the corresponding twistor spinor and
$u:=\left<\phi,\phi\right>$, then the orbifold 
$\overline{M}=(M-\zphi,\overline{g}=u^{-2}g) $ is a Ricci flat
Riemannian orbifold carrying a parallel spinor (recall that $\zphi$ is
the zero set of $\phi$).
This follows from the analogue of 
the Remark in the previous section, applied to orbifolds. 
If $\zphi=\{p_1,\ldots,p_m\}$ then 
$\overline{M}$ has $m$ ends, at any of which the metric $\bar g$ is 
asymptotically locally
Euclidean (i.e. asymptotic to $\R^n/\Gamma_{p_j}$) in
inverted normal coordinates.
As in the Riemannian case we first show that if $m>1$ there is a
geodesic line in $\ol M$. This geodesic line is obtained as a limit of
segments connecting points approaching two different ends (see, for
example, \cite[Ch. 9]{Pe}). 
Let us remark that the corresponding segments can be chosen to avoid the
singularities. Indeed, first note that the singularity set $\Sigma$ is of 
codimension at least 2 (because of the orientability of $M$), so 
$\ol M\smallsetminus\Sigma$ is connected, so we can choose the endpoints of
 the segments to be smooth points. On the other hand, no segment (i.e., 
minimizing geodesic) can touch the singular set without being entirely 
contained in it, otherwise it wouldn't be minimizing 
(the singular set is totally geodesic, for symmetry reasons).
So the singularities play no role in the
process, which carries over exactly as in the smooth case and we get a
geodesic line (which avoids singularities, as well).
On the other hand, there is an orbifold generalization of the
splitting Theorem of Cheeger and Gromoll:

\begin{prop}{\em (\cite[Theorem 1]{BoZh94})}\label{pro:cheegergromoll}
If $\ol M$ is a complete Riemannian orbifold of dimension
$n$ with non-negative Ricci curvature carrying a 
geodesic line, then $\ol M$ is isometric to
the product $\R \times N$ of the real line with
a Riemannian orbifold $N$ of nonnegative 
Ricci curvature.
\end{prop}

\noindent Therefore, if we have more than one end, there should be such a
splitting. But $N\times \R$ is ALE, so the curvature must tend to zero
for large $t$ along the line $\{(x,t)\ |\ t\in\R\}$, for any fixed
$x\in N$. Therefore the curvature of $N$ must be identically zero. On
the other hand, if the ends of that geodesic line are different ends
of $\ol M$, then $N$ must be compact. So $N$ is a flat Riemannian
orbifold.
On the other hand, every end of $\ol M$ corresponds to a removed
(singular) point of singularity group $\Gamma$, therefore $N\simeq
S^{n-1}/\Gamma$, which would imply the existence of a flat metric on
the sphere, contradiction (here we need $n>2$).\\
So there is only one end, hence the zero set of $\phi$ contains only the
point $p$, as claimed.\\
The inequality in the conclusion of the
Theorem is proved using growth estimates for the volume of balls (more
precisely, an extension to orbifolds of Bishop's Theorem):


\begin{prop}{\em (\cite[Proposition 20]{Bo93})}\label{pro:bishop}
Let $\ol M$ be a complete Riemannian orbifold with singular
set $\sing$ and non-negative
Ricci curvature ${\rm Ric} \ge 0.$ 
Then for every point 
$p \in\ol M$ the function
$$r \in (0,\infty) \mapsto \psi (r):= \frac{{\rm vol} B(p,r)}{\omega_n r^n}$$
is non-increasing. 
Here $ {\rm vol} B(p,r)$ is the volume
of the geodesic 
ball $B (p,r)$ of radius $r$ around $p$ and $\omega_n$
is the volume of a unit ball in Euclidean space $\R^n.$
Furthermore 
$$ \lim_{r \to 0} \psi (r)= \frac{1}{\# \Gamma_p}$$
where $\Gamma_p$ is the isotropy subgroup of the
point $p$. Moreover, if, for some $r>0$, $\psi(r)=\frac{1}{\# \Gamma_p}$,
then $B (p,r)$ is isometric to the quotient $B_r^n/\Gamma_p$ of the ball of
radius $r$ in $\R^n$ with $\Gamma_p$.
\end{prop}


\noindent Following the same estimates as in \cite{KR98}, we can show that,
under very weak assumptions, the volume of large balls in an ALE
spoace approaches the one in the flat model:

\begin{lem}\label{lem:volume}
Let $(\ol M,g)$ be an ALE Riemannian orbifold 
(with corresponding group $\Gamma$) of order $(\tau,\mu)$, with 
$\tau>1,\mu\ge 0$.
Then the function $\psi$ (representing the relative volume of a ball in $\ol M$
w.r.t. the Euclidean space) satisfies:
$$\lim_{r\ra\infty}\psi(r)=\frac{1}{\#\Gamma}.$$
\end{lem}


\noindent We apply these inequalities for the proof of Theorem \ref{main}: now,
$\bar M$ is an ALE Ricci-flat orbifold, whose end is asymptotic to
$\R^n/\Gamma$ (Here, $\Gamma$ is actually the group $\Gamma_{p_0}$ 
corresponding to the zero point $p_0$ of the twistor spinor $\phi$).
Choose a point $p \in \ol M$ (if $p$ is a
smooth point, then $\Gamma_p=\{1\}$).
Since the function $\psi (p;r)=\frac{{\rm vol} B(p,r)}{\omega_n r^n}$ is non-increasing it follows
from Lemma \ref{lem:volume} and Proposition \ref{pro:bishop} that 
$$\# \Gamma_p \le \#\Gamma.$$
The equality occurs
if and only if 
$\ol M$ is isometric to $\R^n/\Gamma_p$, therefore $\Gamma\simeq\Gamma_p$ 
and $M$ is conformally equivalent to $S^n/\Gamma$ (where $\Gamma$ acts by 
isometries and has two mutually opposed fixed points). This finishes
the proof of the Theorem.
\end{proof}









\obs We can show that the vanishing locus $\{p_1,\dots,p_k\}$ 
of a twistor spinor on a
compact orbifold contains only singular points without using
Proposition \ref{pro:cheegergromoll}. Indeed, from the volume
estimates given in Proposition~\ref{pro:bishop} we obtain
$$1\ge \frac{1}{\#\Gamma_{p_1}}+\cdots+\frac{1}{\#\Gamma_{p_k}},$$
so if any $p_i$ is smooth, it must be unique and we get equality in
the volume estimate, so the metric is flat. This gives a simpler proof
of the Lichnerowicz' Theorem \cite{Li}.\\

\noindent Note that the above inequality does not exclude the existence of
multiple {\em singular} zeros of a twistor spinor, so the use of
Proposition \ref{pro:cheegergromoll} is necessary.


\begin{cor}\label{nontriv}
Let $(M,g)$ be
a compact Riemannian spin orbifold with finite singularity
set carrying a twistor spinor with non-empty zero set.
Then either $M$ is a quotient of $S^n$ or it is a 
{\em non-trivial} Riemannian orbifold.
\end{cor}

\obs We prefer to use the attribute {\em non-trivial} to denote an
orbifold that can not be obtained
as a quotient of a smooth manifold, as opposed to the slightly
misleading term {\it bad orbifold}, sometimes used in the literature 
to denote the same object.


\begin{proof}
If $M= N/\Gamma$ for a compact Riemannian spin manifold
$N$ then also $N$ carries a twistor spinor with non-empty
zero set. Hence it is a round sphere.
\end{proof}

\noindent Cross-application of the inequality in Theorem \ref{main} leads to:

\begin{cor} If a non-trivial Riemannian orbifold has two (non-trivial) twistor spinors with zero, then they have the same zero set.
\end{cor}

\noindent We can actually say a little more about the sphere quotients that
admit twistor spinors with zeros: 
they have exactly $2$ singular points. Moreover, we can characterize
them (under some technical assumption) merely by the flatness of the conformal
structure:

\begin{prop}\label{confflatex} Any compact, conformally flat orbifold $M$, with isolated singularity set, carrying a twistor 
spinor with zero is conformally equivalent to $S^n/\Gamma$, where $\Gamma\subset SO(n)$ 
acts on $\R^n\subset \R^{n+1}\supset S^n$.
\end{prop}

\begin{proof}
We keep the notations as above. $\ol M$ is a flat, complete orbifold with 
singularity set $\bar\Sigma$. $\ol M\sm\bar\Sigma$ is a manifold, therefore 
it has a universal covering $Q$. Now, in a neighbourhood $U$ of a singular point
$P\in\bar\Sigma$, the map $\pi:Q\ra\ol M\sm\bar\Sigma$ looks like a disjoint
union
of connected coverings of $\R^n\sm\{0\}/\Gamma_P$, where $\Gamma_P$ is the 
singularity group of $P$. Each of these connected components $V_i$ of 
$\pi^{-1}(U)$ can be completed by adding a point $P_i$ to $V_i$. Furthermore, 
$\bar Q:=Q\cup\{P_i\ |\ P\in\bar\Sigma\}$ can be given the structure of a 
flat orbifold, and $\pi:\bar Q\ra\ol M$ is then a ramified covering of
orbifolds.
Since $Q$ is simply connected and flat, any parallel vector field on $Q$ can be
globally extended to a parallel one on $\bar Q$. But $V_i$ is isometric to some
quotient of 
a ball, and its tangent space is spanned by parallel vector fields if and
only if it is a ball itself. So all $V_i$'s are isometric to balls and 
$V_i\cup\{P_i\}$ are smooth balls. So $\bar Q$ is a simply connected smooth
manifold, flat 
and complete, thus $\bar Q\simeq\R^n$. So it has, like $\ol M$, one end, 
therefore the covering $\pi$ is finite (the end of $\ol M$ has finite 
fundamental group).
So $\ol M\simeq\R^n/G$, where $G\subset\mbox{Isom}(R^n)$ is finite and has 
only isolated fixed points. As each $g\in G$ may have at most one fixed point,
the set $\Sigma$ of all fixed points is finite, and
$\bar\Sigma\subset\pi(\Sigma)$.
Let $K$ be the convex envelop of $\Sigma$ and let $P\in\Sigma$ be one of the
vertices of the polyhedron $K$. Then $\Sigma\cap\bar H=\{P\}$ for a certain 
``exterior'' closed halfspace $\bar H$. Let $P'\ne P$ be another point in 
$\Sigma$. Then $P'':=\sum_{k=1}^{N_g}g^k(P')$ is another fixed point for $g$, 
where $g\in G\sm\{1\}$ is such that $g(P)=P$ and $g^{N_g}=1$. But $P''\in
\R^n\sm\bar H$, therefore $P''\ne P$. In this case $g$ would have two distinct
fixed points, so it would leave the whole line connecting them fixed, 
contradiction with the finiteness of $\Sigma$.
So $\Sigma$ consists of only one point and the conclusion easily follows.
\end{proof}

\noindent So the equality case in our Theorem characterizes the conformally flat orbifolds
carrying twistor spinors with zero (which actually turn out to be sphere 
quotients).


\section{Examples}\label{sec:examples}


At our knowledge, no examples of parallel or Killing spinors
(corresponding to non-vanishing solutions of the twistor equation) on
compact, non-smooth, orbifolds are known so far. The following
examples all have exactly one zero and, as pointed out by the
Corollary \ref{nontriv}, they are not only non-smooth, but equally
non-trivial (they are not covered by a manifold).
Let $\Gamma$ be a finite subgroup of the group ${\rm Sp}(1)={\rm SU}(2)$
acting
freely on $\mathbb{H}\sm\{0\},$ here 
$\mathbb{H} \cong \R^4$ is the space of quaternions. Then a
hyperk\"ahler ALE space (of order $(4,\infty)$)
for this group is also called 
{\em gravitational instanton} in the physics literature. 
Outside a compact set the metric is asymptotic to
the quotient $\mathbb{H}/\Gamma$ with order $(4,\infty)$
(cf. Definition~\ref{def:ale}) and in addition also the
hyperk\"ahler structure is asymptotic to the Euclidean
hyperk\"ahler structure, cf. \cite[Definition 7.2.1]{Jo}. 
A hyperk\"ahler space is in particular Ricci flat,
therefore if follows from the rigidity part of
the Bishop-Gromov volume comparison theorem \cite{KR98}
that the group is non-trivial unless the metric is flat.


\begin{exa}\rm\label{exa:eguchi-hanson}
The first examples of a gravitational instanton are the {\em Eguchi-Hanson spaces}
$\left(M_{EH},g_{EH}\right),$
for the subgroup $\Gamma=\Z_2=\{ \pm 1\}.$ These can be given explicitly, 
cf. \cite{EH}, \cite[Example 7.2.2]{Jo}. The space $M_{EH}$ is the blow-up
of $\comp^2 /\Z_2$ at $0,$ it can be identified with the cotangent bundle 
$T^*\mathbb{CP}^1
\cong T^*S^2$
of the $1$-dimensional complex projective space resp. the $2$-sphere.
Hence this space is simply-connected and spin.
It is shown in \cite{KR95} that one can conformally compactify
the Eguchi-Hanson space to a compact orbifold (with 
$C^{\infty}$-metric) with one singular point 
$\pinfty$ whose singularity group is 
$\Z_2.$ The existence of this compactification (at least as
a $C^3$-metric) follows also from
Theorem~\ref{thm:parconf}. \\
Hence we obtain a compact $4$-dimensional Riemannian spin orbifold $(M,g)$ 
with one
singular point $\pinfty$ whose singularity group is $\Z_2$ carrying
two linearly independent twistor spinors $\psi_1,\psi_2.$
The singularity point is the unique zero point of $\psi_1,\psi_2.$ 
\end{exa}


\begin{exa}\rm\label{exa:gibbons-hawking}
Gibbons and Hawking generalized the Eguchi-Hanson construction and
obtained hyperk\"ahler ALE spaces asymptotic to $\HH/\Z_k$ for all
$k \ge 2,$ cf. \cite{GH}. Finally Kronheimer described in
\cite{Kr1} and \cite{Kr2} the construction and classification of 
hyperk\"ahler ALE spaces asymptotic to $\HH/\Gamma$ for nontrivial
finite subgroup $\Gamma \subset {\rm SU}(2),$
see also \cite[Theorem 7.2.3]{Jo}. 
\end{exa}


\begin{exa}\rm\label{exa:calabi}
On $\comp^m\sm\{0\}$ the group $\Z_m$ generated by
complex multiplication with $\zeta=\exp(2\pi \sqrt{-1}/m)$ acts
freely. The blow-up $X$ of $\comp^m/\Z_m$ at $0$  can
be identified with a complex line bundle over the 
$(m-1)$-dimensional complex projective space
$\mathbb{CP}^{m-1}.$
One can explicitly write down an ALE K\"ahler metric
with holonomy ${\rm SU}(m)$, cf. \cite{Ca}, \cite{FG}
or \cite[Example 8.2.5]{Jo}. Using Theorem~\ref{thm:parconf}
or the explicit form given in \cite{KR97} one obtains
a compact spin orbifold 
of real dimension $n=2m$ with one singular point $\pinfty$ and singularity
group $\Z_m.$ \\
This orbifold carries a spin structure with a 
two-dimensional space of twistor spinors whose common zero point
is the singular point.
\end{exa}



\noindent All the examples described so far are defined using a 
{\em crepant resolution} of the isolated singularity $0 \in \comp^2.$
For a definition of a crepant resolution see
\cite[ch. 6.4]{Jo}. In particular for a crepant resolution
the first Chern class vanishes, hence
one can show that a resolution of
$\comp^m /\Gamma$ carries a Ricci flat ALE K\"ahler metric only if
the resolution is crepant, cf. \cite[Proposition 8.2.1]{Jo}. 
In complex dimension 
$2$ there is a unique crepant resolution of $\comp^m/\Gamma$ (here 
$\Gamma\subset SU(n)$ and acts freely on the sphere in $\C^n$), 
and in dimension $3$, for each $\Gamma$ as above, 
there exist crepant resolutions (not unique, in general), but 
in higher dimensions this does no longer hold, for example
$\comp^4/\{\pm1\}$ does not carry any crepant resolution, 
cf. \cite[Example 6.4.5]{Jo}.
For crepant resolutions $M$ of the quotient singularity $\comp^m/\Gamma$, 
if there is an ALE K\"ahler metric on $M$, there is a solution of the 
ALE Calabi conjecture, as stated and proved by D. Joyce, \cite[Theorem
8.2.3]{Jo}, 
\cite[Theorem 8.2.4]{Jo}: In every 
K\"ahler class of ALE K\"ahler metrics (with order $(2m,\infty)$)
there is a unique Ricci-flat K\"ahler metric with holonomy 
${\rm SU}(m).$ Then we can conclude from \cite[Corollary 3.6.3]{Jo}
that this metric is spin and carries a $2$-dimensional space
of parallel spinors. The resolution is simply-connected, hence
there is only one spin structure. Hence starting from a
crepant resolution of the quotient sigularity
$\comp^m/\Gamma$ with an ALE K\"ahler metric 
(of order $(2m,\infty)$) one obtains
a compact spin orbifold of dimension $n=2m$
with a $C^{2m-1}$-Riemannian
metric which carries a $2$-dimensional
space of twistor spinors. All twistor spinors have exactly one
zero in the only singularity point.
In general, 
the existence of ALE K\"ahler metrics on crepant resolutions is an
open question, however, in dimension 3, ALE K\"ahler (non-Ricci flat)
metrics have been constructed on crepant resolutions of $\C^3/\Gamma$
using symplectic quotient techniques, \cite{SI},\cite{deg}.
Therefore we can add the following class of examples of 6-dimensional
orbifolds:


\begin{exa}\rm\label{exa:dim3}
Let $\Gamma\subset SU(3)$ be any finite subgroup that
acts freely on $S^5\subset \C^3$. Then there exists a crepant resolution $M$ of $\C^3/\Gamma$, admitting an ALE K\"ahler metric
of order $(6,\infty)$, \cite{deg},
and, therefore, an ALE Ricci-flat one as well, cf. \cite{Jo}. 
The one-point compactification of $M$ is then spin, has one singular point 
at $\infty$ (of singularity group $\Gamma$), and admits two linearly independent
twistor spinors that vanish at 
that singular point.
\end{exa}


\obs Theorem \ref{thm:parconf} shows only that the metric has some finite
regularity $C^k$ around the zero of the twistor spinor, depending on
the decay rate of the ALE K\"ahler Ricci-flat metric on the
complement. 
It is unknown whether this can be improved, possibly after changing
the metric in the conformal class.


























\end{document}